\documentstyle[12pt]{article}

\begin{document}

\baselineskip16.5pt

\newtheorem{definition}{Definition $\!\!$}[section]
\newtheorem{prop}[definition]{Proposition $\!\!$}
\newtheorem{lem}[definition]{Lemma $\!\!$}
\newtheorem{corollary}[definition]{Corollary $\!\!$}
\newtheorem{theorem}[definition]{Theorem $\!\!$}
\newtheorem{example}[definition]{\it Example $\!\!$}
\newtheorem{remark}[definition]{Remark $\!\!$}

\newcommand{\nc}[2]{\newcommand{#1}{#2}}
\newcommand{\rnc}[2]{\renewcommand{#1}{#2}}

\def\inbar{\,\vrule height1.5ex width.4pt depth0pt}
\def\IN{\relax\,\hbox{$\inbar\kern-.3em{\rm N}$}}
\def\IZ{\relax\,\hbox{$\inbar\kern-.3em{\rm Z}$}}
\def\IC{\relax\,\hbox{$\inbar\kern-.3em{\rm C}$}}
\def\IR{\relax\,\hbox{$\inbar\kern-.3em{\rm R}$}}

\nc{\bpr}{\begin{prop}}
\nc{\bth}{\begin{theorem}}
\nc{\ble}{\begin{lem}}
\nc{\bco}{\begin{corollary}}
\nc{\bre}{\begin{remark}}
\nc{\bex}{\begin{example}}
\nc{\bde}{\begin{definition}}
\nc{\ede}{\end{definition}}
\nc{\epr}{\end{prop}}
\nc{\ethe}{\end{theorem}}
\nc{\ele}{\end{lem}}
\nc{\eco}{\end{corollary}}
\nc{\ere}{\hfill\mbox{$\Diamond$}\end{remark}}
\nc{\eex}{\end{example}}
\nc{\epf}{\hfill\mbox{$\Box$}}
\nc{\ot}{\otimes}
\nc{\bsb}{\begin{Sb}}
\nc{\esb}{\end{Sb}}
\nc{\ct}{\mbox{${\cal T}$}}
\nc{\ctb}{\mbox{${\cal T}\sb B$}}
\nc{\bcd}{\[\begin{CD}}
\nc{\ecd}{\end{CD}\]}
\nc{\ba}{\begin{array}}
\nc{\ea}{\end{array}}
\nc{\bea}{\begin{eqnarray}}
\nc{\eea}{\end{eqnarray}}
\nc{\be}{\begin{enumerate}}
\nc{\ee}{\end{enumerate}}
\nc{\beq}{\begin{equation}}
\nc{\eeq}{\end{equation}}
\nc{\bi}{\begin{itemize}}
\nc{\ei}{\end{itemize}}
\nc{\kr}{\mbox{Ker}}
\nc{\te}{\!\ot\!}
\nc{\pf}{\mbox{$P\!\sb F$}}
\nc{\pn}{\mbox{$P\!\sb\nu$}}
\nc{\bmlp}{\mbox{\boldmath$\left(\right.$}}
\nc{\bmrp}{\mbox{\boldmath$\left.\right)$}}
\rnc{\phi}{\mbox{$\varphi$}}
\nc{\LAblp}{\mbox{\LARGE\boldmath$($}}
\nc{\LAbrp}{\mbox{\LARGE\boldmath$)$}}
\nc{\Lblp}{\mbox{\Large\boldmath$($}}
\nc{\Lbrp}{\mbox{\Large\boldmath$)$}}
\nc{\lblp}{\mbox{\large\boldmath$($}}
\nc{\lbrp}{\mbox{\large\boldmath$)$}}
\nc{\blp}{\mbox{\boldmath$($}}
\nc{\brp}{\mbox{\boldmath$)$}}
\nc{\LAlp}{\mbox{\LARGE $($}}
\nc{\LArp}{\mbox{\LARGE $)$}}
\nc{\Llp}{\mbox{\Large $($}}
\nc{\Lrp}{\mbox{\Large $)$}}
\nc{\llp}{\mbox{\large $($}}
\nc{\lrp}{\mbox{\large $)$}}
\nc{\lbc}{\mbox{\Large\boldmath$,$}}
\nc{\lc}{\mbox{\Large$,$}}
\nc{\Lall}{\mbox{\Large$\forall$}}
\nc{\bc}{\mbox{\boldmath$,$}}
\rnc{\epsilon}{\varepsilon}
\rnc{\ker}{\mbox{\em Ker}}
\nc{\ra}{\rightarrow}
\nc{\ci}{\circ}
\nc{\cc}{\!\ci\!}
\nc{\T}{\mbox{\sf T}}
\nc{\can}{\mbox{\em\sf T}\!\sb R}
\nc{\cnl}{$\mbox{\sf T}\!\sb R$}
\nc{\lra}{\longrightarrow}
\nc{\M}{\mbox{Map}}
\rnc{\to}{\mapsto}
\nc{\imp}{\Rightarrow}
\rnc{\iff}{\Leftrightarrow}
\nc{\bmq}{\cite{bmq}}
\nc{\ob}{\mbox{$\Omega\sp{1}\! (\! B)$}}
\nc{\op}{\mbox{$\Omega\sp{1}\! (\! P)$}}
\nc{\oa}{\mbox{$\Omega\sp{1}\! (\! A)$}}
\nc{\inc}{\mbox{$\,\subseteq\;$}}
\nc{\de}{\mbox{$\hD $}}
\nc{\spp}{\mbox{${\cal S}{\cal P}A({\IC}^2_p)\# A(GL_{q,p}(2))$}}
\nc{\dr}{\mbox{$\hD _{R}$}}
\nc{\dsr}{\mbox{$\hD _{\cal R}$}}
\nc{\m}{\mbox{m}}
\nc{\0}{\sb{(0)}}
\nc{\1}{\sb{(1)}}
\nc{\2}{\sb{(2)}}
\nc{\3}{\sb{(3)}}
\nc{\4}{\sb{(4)}}
\nc{\5}{\sb{(5)}}
\nc{\6}{\sb{(6)}}
\nc{\7}{\sb{(7)}}
\nc{\hsp}{\hspace*}
\nc{\nin}{\mbox{$n\in\{ 0\}\!\cup\!{\IN}$}}
\nc{\al}{\mbox{$\alpha$}}
\nc{\bet}{\mbox{$\beta$}}
\nc{\ha}{\mbox{$\alpha$}}
\nc{\hb}{\mbox{$\beta$}}
\nc{\hg}{\mbox{$\gamma$}}
\nc{\hd}{\mbox{$\delta$}}
\nc{\he}{\mbox{$\varepsilon$}}
\nc{\hz}{\mbox{$\zeta$}}
\nc{\hs}{\mbox{$\sigma$}}
\nc{\hk}{\mbox{$\kappa$}}
\nc{\hm}{\mbox{$\mu$}}
\nc{\hn}{\mbox{$\nu$}}
\nc{\la}{\mbox{$\lambda$}}
\nc{\hG}{\mbox{$\Gamma$}}
\nc{\hD}{\mbox{$\Delta$}}
\nc{\th}{\mbox{$\theta$}}
\nc{\Th}{\mbox{$\Theta$}}
\nc{\ho}{\mbox{$\omega$}}
\nc{\hO}{\mbox{$\Omega$}}
\nc{\hp}{\mbox{$\pi$}}
\nc{\hP}{\mbox{$\Pi$}}
\nc{\bpf}{{\it Proof.~~}}
\nc{\slq}{\mbox{$A(SL\sb q(2))$}}
\nc{\fr}{\mbox{$Fr\llp A(SL(2,\IC))\lrp$}}
\nc{\slc}{\mbox{$A(SL(2,\IC))$}}
\nc{\af}{\mbox{$A(F)$}}
\nc{\w}{\mbox{\sf w}}
\nc{\r}{\mbox{\sf r}}
\nc{\s}{\mbox{\sf s}}
\nc{\fu}{\mbox{\sf u}}
\nc{\fv}{\mbox{\sf v}}

\nc{\fa}{\mbox{\sf a}}
\nc{\fb}{\mbox{\sf b}}
\nc{\fc}{\mbox{\sf c}}
\nc{\fd}{\mbox{\sf d}}
\nc{\D}{\mbox{\sf D}}
\nc{\x}{\mbox{\sf x}}
\nc{\y}{\mbox{\sf y}}

\def\esl{{\mbox{$E\sb{\frak s\frak l (2,{\IC})}$}}}
\def\esu{{\mbox{$E\sb{\frak s\frak u(2)}$}}}
\def\zf{{\mbox{${\IZ}\sb 4$}}}
\def\zt{{\mbox{$2{\IZ}\sb 2$}}}
\def\ox{{\mbox{$\Omega\sp 1\sb{\frak M}X$}}}
\def\oxh{{\mbox{$\Omega\sp 1\sb{\frak M-hor}X$}}}
\def\oxs{{\mbox{$\Omega\sp 1\sb{\frak M-shor}X$}}}
\def\st{\stackrel}
\def\Fr{\mbox{Fr}}
\def\gal{-Galois extension}
\def\hge{Hopf-Galois extension}

\begin{center}
\vspace*{1.0cm}

{\LARGE{\bf  
Frame, cotangent and tangent bundles of the quantum plane}}\footnote{Modified
version of the talk given by R.M. at the workshop ``Lie Theory and its 
Applications in Physics II", Clausthal, August 1997.} 

\vskip 1.5cm

{\large {\bf Piotr M.~Hajac}}\footnote{
Partially supported by the NATO Postdoctoral Fellowship and by 
the KBN grant \mbox{2 P301 020 07}. On~leave from:
Department of Mathematical Methods in Physics, 
Warsaw University, ul.~Ho\.{z}a 74, Warsaw, \mbox{00--682~Poland}.
 http://info.fuw.edu.pl/KMMF/ludzie$\underline{~~}$ang.html 
e-mail: pmh@fuw.edu.pl} 

\vskip 0.5 cm 

Department of Applied Mathematics and Theoretical Physics\\
University of Cambridge, Silver Street, Cambridge CB3 9EW, England\\
http://www.damtp.cam.ac.uk/user/pmh33\\ 
e-mail: pmh33@amtp.cam.ac.uk \\

\vskip 1cm

{\large{\bf Rainer Matthes}}

\vskip 0.5cm

Institut f\"ur Theoretische Physik der Universit\"at Leipzig\\
Augustuplatz 10/11, D-04109 Leipzig, Germany\\
e-mail: matthes@itp.uni-leipzig.de

\end{center}

\vspace{1 cm}

\begin{abstract}
We construct a quantum frame bundle of the quantum plane ${\IC}
^2_p$ by requiring that a $GL_{q,p}(2)$-covariant differential calculus
 on ${\IC}^2_p$ be isomorphic as a bimodule to the space
of sections of the associated quantum cotangent bundle. 
We also construct the section space of the associated quantum
tangent bundle, and show that it is naturally dual to the differential
calculus.
\end{abstract}

\vspace{1 cm}

\section{Introduction}

The notion of a \hge\ replaces the classical concept of a principal bundle
the same way Hopf algebras replace groups. This is why we think of objects dual 
to \hge s as quantum principal bundles. An obvious next step is to consider
associated quantum vector bundles. More precisely, what we need here is a
replacement of the module of sections of an associated vector bundle. In the
classical case such sections can be equivalently described  as ``functions
of type $\varrho$" from the total space of a principal bundle to a vector 
space. We follow this construction in the quantum case by considering bimodules
associated with \hge s (see~\cite[Appendix~B]{d1} and~\cite{d2}).

The aim of this note is to present an example of a \hge\ which can be considered 
the dual of a $GL_{q,p}(2)$-frame bundle of the quantum plane ${\IC}^2_p$.
It is constructed in such a way that the naturally associated cotangent bimodule 
is isomorphic to one of the known left covariant differential calculi on 
${\IC}^2_p$. The main claim of this work is contained in the last section, 
where we determine the left and right bimodule duals of the aforementioned
differential calculus, and prove that they are both isomorphic to the associated 
tangent bimodule.

We assume that the reader is familiar with the Hopf-algebra terminology
and Sweedler notation $\Delta(f)=\sum f_1\ot f_2$, etc.

\section{Preliminaries}

We recall here several definitions and known results to be used in the sequel.

\bde
Let $H$ be a Hopf algebra, let $P$ be a right $H$-comodule algebra with a
structure map $\Delta_R:P\ra P\ot H$, and let $\iota:B\ra P$ be an algebra
embedding such that $\iota(B)=P^{coH}:=\{f\in P|\Delta_R(f)=f\ot 1\}$. 
If the canonical map 
\[
P\ot_{B}P\lra P\ot H,\;\;\; f\ot f'\lra \sum ff'_0\ot f'_1\, ,
\]
 is bijective, $B\inc P$ is called an $H$-Galois extension.
\ede
\bde
An $H$-Galois extension $B\inc P$ is called cleft if there exists a linear
map (called a cleaving map) $j:H\ra P$ which is invertible in the convolution 
algebra $Hom(H,P)$ and fulfills 
$\Delta_R\ci j=(j\ot id)\ci \Delta$ ($H$-colinearity).
\ede
Cleft extensions are equivalent to crossed products $B\#_{(\triangleright,\sigma)}
H$. The latter are determined by a pair $(\triangleright, \sigma)$ 
(left cocycle action) defined as follows:
\bde
A left cocycle action of $H$ on $B$ is a pair $(\triangleright,\sigma)$, 
where $H\otimes B\st{\triangleright}{\lra} B$ and $H\ot H\st{\sigma}{\lra} B$ 
are linear maps with the properties\\
(i) $ h\triangleright(ab)=\sum (h_1\triangleright a)(h_2\triangleright b)$,\\
(ii) $h\triangleright 1=\varepsilon(h)1$,\\
(iii) $1\triangleright b=b$,\\
(iv) $\sigma(h\ot 1)=\sigma(1\ot h)=\varepsilon(h)1$,\\
(v) $\sum h_1\triangleright\sigma(h'_1\ot h''_1)\sigma(h_2\ot h'_2h''_2)=
\sum \sigma(h_1\ot h'_1)\sigma(h_2h'_2\ot h'')$,\\
(vi) $\sum(h_1\triangleright(h'_1\triangleright b))\sigma(h_2\ot h'_2)=
\sum\sigma(h_1\ot h'_1)((h_2h'_2)\triangleright b)$,\\
(vii) $\sigma$ is convolution invertible.
\ede
\bpr
Let $(\triangleright,\sigma)$ be a left cocycle action. The bilinear map defined by
\[(b\ot h)(b'\ot h'):=\sum b(h_1\triangleright b')\sigma(h_2\ot h'_1)\ot h_3h'_2\]
makes $B\ot H$ a unital associative algebra.
\epr
This algebra, denoted by $B\#_{(\triangleright,\sigma)}H$,
 is called a crossed product
of $B$ and $H$. $B$ is embedded as an algebra in $B\#_{(\triangleright,\sigma)}H$ 
via $id\ot 1$. The homomorphism $id\ot\Delta$ makes the latter a right $H$-comodule 
algebra. $B$ is the corresponding subalgebra of coinvariants, and $1\ot id$ is a
cleaving map. Conversely, for a cleft $H$\gal\  $B\inc P$ with a
cleaving map $j:H\ra P$, one defines a left cocycle action
$(\triangleright,\sigma)$ by
\[
h\triangleright b=\sum j(h_1)bj^{-1}(h_2),\;\;\;
\sigma(h\ot h')=\sum j(h_1)j(h'_1)j^{-1}(h_2h'_2),
\]
where $j^{-1}$ is the convolution inverse of $j$. It turns out that
$\phi:B\#_{(\triangleright,\sigma)} H\ra P$, $b\ot h\to bj(h)$, is an isomorphism 
of right $H$-comodule algebras.
Smash products are defined as crossed products with the trivial cocycle $\sigma$, 
i.e., $\sigma=\varepsilon\ot\varepsilon ~1$. In this case, the above conditions
(iv), (v) and (vii) become trivial, whereas (vi) becomes 
\[
(vi')\;\;\; h\triangleright(h'\triangleright b)=(hh')\triangleright b.
\]
The map $\triangleright$ with the properties (i), (ii), (iii) and (vi') is called 
a left action of $H$ on $B$. The smash product of $B$ and $H$ is usually denoted by
$B\# H$, though it may still depend on the choice of the left action.
(For information about these notions and facts see 
\cite{bcm,blmo,dota0,m-s,maj23.5,mo,s3} and references therein.)

For an $H$\gal\ $B\inc P$ and a right $H$-comodule $F$ ($\rho_R:F\ra F\ot H$),
we think of the associated bimodule of colinear maps (intertwiners) 
$_{\rho_R}Hom_{\Delta_R}(F,P):=\{\ell\in Hom(F,P)\,|
\,\Delta_R\ci\ell=(\ell\ot id)\ci\rho_R\}$ as the $B$-bimodule of sections
of the associated quantum vector bundle. The bimodule structure of 
$_{\rho_R}Hom_{\Delta_R}(F,P)$ is defined by
$(b\ell)(f)=b(\ell(f)),~~(\ell b)(f)=(\ell(f))b$.
Note that the classical cotangent bundle corresponds to the data 
$F={{\IR}^n}$, $H=C^{\infty}(GL(n))$, $\rho_{R}(e_i)=e_j\ot {g^j}_i$, where
$\{e_i\}$ is a basis of ${\IR}^n$, and the matrix elements ${g^j}_i$
are considered functions on $GL(n)$.
The corresponding data for the tangent bundle is $F={{\IR}^n}^*,~~H=C^{\infty}
(GL(n)), ~~\tilde{\rho}_R(e^i)=e^j\ot S({g^i}_j)$. Here $\{e^i\}$ is a basis of 
${{\IR}^n}^*$, and $S$ is the antipode in $C^{\infty}(GL(n))$. This motivates
our constructions and terminology in the subsequent sections.

Let us now recall the definition of the coordinate ring of the quantum plane 
${\IC}^2_q$ and two-parameter quantum general linear group $GL_{q,p}(2)$.
$A({\IC}^2_p)$ is defined as the unital associative algebra over ${\IC}$ 
generated by $x,y$ with relations
\beq\label{cp2} 
xy=pyx,\hspace{.5cm} p\in{\IC}\setminus\{0\}.
\eeq
$A(GL_{q,p}(2))$ is defined as the unital associative algebra over $\IC$
generated by $a,b,c,d,D^{-1}$ with relations
\beq 
ab=qba,\hspace{.2cm} ac=pca,\hspace{.2cm} bd=pdb,
\hspace{.2cm} cd=qdc,\hspace{.2cm} bc=\frac{p}{q}cb,\hspace{.2cm}ad=da+
(q-p^{-1})bc\label{rel1}
\eeq
\beq\label{rel3} 
(ad-qbc)D^{-1}=D^{-1}(ad-qbc)=1 ,
\eeq
where $p,q\in{\IC}\setminus\{0\}$.
(Notice that, compared with \cite{man4}, we have used $p^{-1},q^{-1}$ instead of 
$q,p$.)
For the inverse of the quantum determinant $D=ad-qbc$
one derives the commutation rules
\beq\label{cot} 
aD^{-1}=D^{-1}a, ~dD^{-1}=D^{-1}d,~ bD^{-1}=qp^{-1}D^{-1}d, ~cD^{-1}
=q^{-1}pD^{-1}c.
\eeq
The Hopf algebra structure of $A(GL_{q,p}(2))$ is defined in terms of the matrix
\[
T=\left(\ba{cc}a & b \\ c & d \ea\right)
\]
of generators by
\beq
\hD(T)=T\stackrel{.}{\ot} T, ~\hD(D^{-1})=D^{-1}\ot D^{-1} ,\label{delt}
\eeq
\beq
\epsilon(T)=I , ~\epsilon(D^{-1})=1,\label{et} 
\eeq
\beq 
S(T)=\left(\ba{cc}d & -q^{-1}b\\
-qc & a\ea\right)D^{-1}, ~S(D^{-1})=D. \label{St} 
\eeq
Here all equations involving $T$ are to be understood componentwise.
$A({\IC}^2_p)$ is a left $A(GL_{q,p}(2))$  and a right 
$A(GL_{p,q}(2))$-comodule algebra with coactions $\hd_L$ and $\hd_R$ 
given on generators by $\hd_L(X)=T\stackrel{.}{\ot} X$ and 
$\hd_R(X^t)=X^t\stackrel{.}{\ot} T \label{col}$, where
\[
X=\left(\ba{c}x\\y\ea\right)
\]
 and $X^t$ is the transposed matrix.

Next, recall from \cite{man4} (cf.~\cite{bdr}) that the $A({\IC}^2_p)$-bimodule 
$\Gamma(T\sp*_q{\IC}^2_p)$ 
generated by
 $\xi,\eta$ with relations
\beq\label{x11}
x\xi=pq\xi x,~~x\eta = (pq-1)\xi y + p\eta x,
\eeq
\beq\label{x22}
y\xi=q\xi y,~~y\eta=pq\eta y
\eeq
is a left $A(GL_{q,p}(2))$-covariant $A({\IC}^2_p)$-bimodule. The 
left coaction
\[
\hat{\hD}_L:\Gamma(T\sp*_q{\IC}^2_p)\lra A(GL_{q,p}(2))\ot 
\Gamma(T\sp*_q{\IC}^2_p)
\]
is given by
\beq
\hat{\hD}_L\left(\ba{c}\xi\\ \eta\ea\right)=\left(\ba{cc}a&b\\c&d\ea\right)
\stackrel{.}{\ot}
\left(\ba{c}\xi\\\eta\ea\right).
\label{dhat}
\eeq
Moreover, with $dx:=\xi$, $dy:=\eta$, $(\Gamma(T\sp*_q{\IC}^2_p),d)$ is a 
first order differential calculus over ${\IC}^2_p$.
(Comparing with \cite{man4}, we have used again $p^{-1},q^{-1}$
instead of $q,p$.)

\section{Frame and cotangent bundles of \boldmath ${\IC}^2_p$}

We are now going to construct a \hge\ playing the role of a quantum frame
bundle of ${\IC}^2_p$. This construction essentially resembles that of
\cite[Example~10.2.12]{m-s}. First we need:
\ble\label{bml}
There exists a well-defined left action of $A(GL_{q,p}(2))$ on $A({\IC}^2_p)$ 
given by the formulas
\beq\label{wx} 
a\triangleright x=(pq)^{-1}x, ~b\triangleright x=0, ~c\triangleright x=((pq)^{-1}-1)y, ~d\triangleright x=p^{-1}x, ~D^{-1}\triangleright x=p^2qx,
\eeq
\beq\label{wy} 
a\triangleright y=q^{-1}y, ~b\triangleright y=0, ~c\triangleright y=0, ~d\triangleright y=(pq)^{-1}y, ~D^{-1}\triangleright y=pq^2y.
\eeq
\ele
\bpf
A direct verification that (\ref{wx})--(\ref{wy}) respect the ideals defining
$A({\IC}^2_p)$ and\linebreak $A(GL_{q,p}(2))$.
\epf\\ \ \\
The above action defines a smash product whose cross
relations reproduce the bimodule structure of $\Gamma(T\sp*_q{\IC}^2)$.
We denote the smash product $A({\IC}^2_p)\# A(GL_{q,p}(2))$ 
corresponding to this action by $A(F_q{\IC}^2_p)$, and call it the coordinate
ring of a quantum frame bundle of the quantum plane ${\IC}^2_p$.  
We have injective algebra homomorphisms 
$id\ot 1:A({\IC}^2_p)\ra A(F_q{\IC}^2_p)$ and 
$1\ot id: A(GL_{q,p}(2))\ra A(F_q{\IC}^2_p)$, and the right 
coaction $\Delta_R=id\ot\Delta$ making $A({\IC}^2_p)$  the algebra 
of coinvariants. The images $x\ot 1,y\ot 1, 1\ot a,1\ot b,\ldots,1\ot D^{-1}$ 
of the generators of $A({\IC}^2_p)$ and $A(GL_{q,p}(2))$ under these embeddings 
are denoted by $\x, \y, \fa,\ldots, \D^{-1}$ respectively.
They are the generators of $A(F_q{\IC}^2_p)$. They  satisfy the
relations (\ref{cp2}), (\ref{rel1}), (\ref{rel3}) and additional cross 
relations coming from the definition of the smash product:
\beq\label{x1}
\x\fa=pq\fa\x,~~\x\fb=pq\fb,~~
\x\fc=(pq-1)\fa\y+p\fc\x,~~\x\fd=(pq-1)\fb\y+p\fd\x,
\eeq
\beq\label{x2}
\y\fa=q\fa\y,~~\y\fb=q\fb\y,~~\y\fc=pq\fc\y,~~\y\fd=pq\fd\y.
\eeq
Furthermore, one can compute
\beq\label{yt}
\x\D^{-1}=p^{-2}q^{-1}\D^{-1}\x,
\hspace*{.3cm}\y\D^{-1}=p^{-1}q^{-2}\D^{-1}\y.
\eeq
As can be shown using Proposition~4.5 of~\cite{maj23.5}, 
$A(F_q{\IC}^2_p)$ is not isomorphic an an algebra to the tensor product 
$A({\IC}^2_p)\ot A(GL_{q,p}(2))$.

Our goal is now to show that the differential calculus defined by (\ref{x11})
and (\ref{x22}), 
which is one of the two left $GL_{q,p}(2)$-covariant first order 
differential calculi on ${\IC}^2_p$ described in \cite{man4}, is isomorphic 
as an $A({\IC}^2_p)$-bimodule to the cotangent bimodule
associated to $A(F_q{\IC}^2_p)$. 

Let $\{e=(1,0),f=(0,1)\}$ be the canonical basis of ${\IC}^2$, and 
 $\rho_R$ be the right corepresentation of $A(GL_{q,p}(2))$ on ${\IC}^2$ 
given by the formula
\beq
\rho_R(e~~f)=(e~~f)\stackrel{.}{\ot}\left(\ba{cc}a&b\\c&d\ea\right).\label{ro}
\eeq
Following arguments of the previous section, we treat the bimodule of colinear
maps\linebreak $_{\rho_R}Hom_{\hD_R}({\IC}^2,A(F_q{\IC}^2_p))$ as the
space of sections of the quantum cotangent bundle of~${\IC}^2_p$.
The main point of this section is contained in the following:
\bpr\label{biso}
$_{\rho_R}Hom_{\hD_R}({\IC}^2,A(F_q{\IC}^2_p))$ and 
$\Gamma(T\sp*_q{\IC}^2_p)$ defined by (\ref{x11})--(\ref{x22}) are isomorphic
as $A({\IC}^2_p)$-bimodules.
\epr
\bpf
First we need to state the following reformulation of
\cite[Proposition~A.7]{bm}:
\ble\label{isol}
Let $B\inc P$ be a cleft $H$\gal\ with a cleaving map $j$, and let
$\rho_R:F\ra F\ot H$ be a right corepresentation of $H$ on $F$.
Then $_{\rho_R}Hom_{\Delta_R}(F,P)$ is isomorphic as a left $B$-module  to
the free module $Hom(F,B)$.
\ele
The isomorphism $\Psi:Hom(F,B)\,\ra\; _{\rho_R}Hom_{\hD_R}(F,P)$ used in the
proof of the above lemma is given by the formulas~\cite{bm}:
\[
(\Psi(u))(\la)=\sum u(\la_0)j(\la_1),\;\;\;
\rho_R(\la)=\sum\la_0\ot\la_1\, .
\]
Denote by $\{\hs_x,\hs_y\}$ the basis of the $A({\IC}^2_p)$-module
$Hom({\IC}^2,A({\IC}^2_p))$ dual to the canonical basis $\{e,f\}$ of
${\IC}^2$. (To simplify notation, here and further on, we identify
$A({\IC}^2_p)$ with its image in $A(F_q{\IC}^2_p)$.)
By Lemma~\ref{isol}, 
$\{\Psi(\hs_x),\Psi(\hs_y)\}$ is a basis of the left
$A({\IC}^2_p)$-module $_{\rho_R}Hom_{\hD_R}({\IC}^2,A(F_q{\IC}^2_p))$.
To fix $\Psi$, we choose a cleaving map $j$ to be $1\ot id$.
One can now check directly that the mapping determined by
\[
\xi\to\Psi(\hs_x),~\eta\to\Psi(\hs_y)
\]
yields the desired isomorphism of bimodules.
\epf

\section{Tangent construction for \boldmath ${\IC}^2_p$}

Following the construction of the associated cotangent bimodule 
$\Gamma(T\sp*_q{\IC}^2_p)$, we define the tangent bimodule:
\bde
Let $\tilde{\rho}_R:{\IC}^2\ra{\IC}^2\ot A(GL_{q,p}(2))$ be the right 
corepresentation given by
\beq\label{roti}
\tilde{\rho}_R(e,f)=(e,f)\stackrel{.}{\ot} S\left(\ba{cc}a&b\\c&d\ea\right)^t, 
\eeq
 where $\{e,f\}$ is the canonical basis of ${\IC}^2$ as above. The bimodule
\[
\Gamma(T_q{\IC}^2_p):=\, 
_{\tilde{\rho}_R}Hom_{\hD_R}({\IC}^2,A(F_q{\IC}^2_p))
\] 
is called the space of sections of the quantum tangent bundle of ${\IC}^2_p$.
\ede
Taking again advantage of Lemma~\ref{isol}, we find a left 
$A({\IC}^2_p)$-module basis of $\Gamma(T_q{\IC}^2_p)$ from the formulas:
\[
\partial_x(\la)=\sum\sigma_x(\la_0)j(\la_1),\;\;\;
\partial_y(\la)=\sum\sigma_y(\la_0)j(\la_1)\, .
\]
Here $\tilde\rho_R(\la)=\sum\la_0\ot\la_1$, and $\sigma_x$, $\sigma_y$ and $j$
are as in the previous section. Moreover, one can directly compute the following:
\ble
The bimodule structure of $\Gamma(T_q{\IC}^2_p)$ is determined by the 
equations:

\beq\label{xdx}
x\partial_x=(q^{-1}p^{-1}-1)\partial_yy+p^{-1}q^{-1}\partial_xx,
\eeq
\beq\label{xdy}
x\partial_y=p^{-1}\partial_yx ,
\eeq
\beq\label{ydx}
y\partial_x=q^{-1}\partial_xy ,
\eeq
\beq\label{ydy}
y\partial_y=p^{-1}q^{-1}\partial_yy .
\eeq
\ele\ \\ 
Caution: Symbols like $\partial_xx$ do not signify differentiation but 
multiplication of $\partial_x$ by $x$ on the right, i.e., 
$(\partial_xx)(\la):=(\partial_x(\la))x$.
\bco
$\{\partial_x,\partial_y\}$ is a basis of the right $A({\IC}^2_p)$-module 
$\Gamma(T_q{\IC}^2_p)$.
\eco

We will now show that there is a natural duality between the 
$A({\IC}^2_p)$-bimodules $\Gamma(T\sp*_q{\IC}^2_p)$ and 
$\Gamma(T_q{\IC}^2_p)$. First, let us recall the definition of
bimodule duals~(cf. \cite{bou}). Let $B$ be an algebra, and let $M$ be a $B$-bimodule. 
Denote by $^*M$ the space all of left $B$-module maps $M\ra B$, and by
$M^*$ the space all of right $B$-module maps $M\ra B$. 
$^*M$ becomes a $B$-bimodule via
\[
(bX)(m)=X(mb),~~~(Xb)(m)=X(m)b,~~~X\in {^*M},~~~ m\in M,~~~ b\in B.
\]
Similarly, the bimodule structure of $M^*$ is given by:
\[
(bX)(m)=bX(m),~~~(Xb)(m)=X(bm).
\]
We call the bimodules $^*M$ and $M^*$ the left and right dual of $M$ respectively. These notions of bimodule duals were used in \cite{bo} to define
Cartan pairs, and thus obtain a general concept of vector field in
Noncommutative Geometry.
Now we can state the main claim of this paper:
\bpr\label{main}
The $A({\IC}^2_p)$-bimodule $\Gamma(T_q{\IC}^2_p)$ is naturally isomorphic to
the left dual of the $A({\IC}^2_p)$-bimodule $\Gamma(T\sp*_q{\IC}^2_p)$. 
 The left and right dual of 
$\Gamma(T\sp*_q{\IC}^2_p)$ are also isomorphic as bimodules.
\epr
\bpf
The first assertion follows from the fact that  the left dual basis of the
basis $\{\xi,\eta\}$ of $\Gamma(T\sp*_q{\IC}^2_p)$   
satisfies (\ref{xdx})-(\ref{ydy}). The dual basis $\{\partial_x^R
,\partial_y^R\}$ of the right 
dual of $\Gamma(T\sp*_q{\IC}^2_p)$  satisfies (\ref{xdy})-(\ref{ydy}) 
and 
\[
x\partial_x^R=(pq)^{-1}((pq)^{-1}-1)\partial_y^Ry+(pq)^{-1}\partial_x^Rx.
\]
Hence $\partial_x\to\partial_x^R,~~\partial_y\to (pq)^{-1}\partial_y^R$ 
defines a bimodule isomorphism proving the second assertion.
\epf


{\it Acknowledgments:}~
It is a pleasure to thank Mico Durdevic, Istvan Heckenberger
 and Axel Sch\"{u}ller for
helpful doscussions. P.M.H.\ is also very grateful for the kind
hospitality of the University of Leipzig.

{\em Erratum:}~
 This paper replaces its earlier version titled ``A Frame Bundle of
the Quantum Plane", preprint NTZ 37/1997. In the previous version, the
construction of the smash
product comodule algebra $A(F_q{\IC}^2_p)$ was flawed. We used $GL_{p,q}(2)$
instead of $GL_{q,p}(2)$. This innocent looking interchange causes the collapse
of the algebra  $A(F_q{\IC}^2_p)$. In particular, $A({\IC}^2_p)$ is no 
longer a subalgebra of $A(F_q{\IC}^2_p)$.

\end{document}